\documentclass[12pt]{article}

\usepackage{amsmath}
\usepackage{amssymb}
\usepackage{amsthm}
\usepackage{mathtools}
\usepackage{indentfirst}

\usepackage{hyperref}

\newtheorem*{lemma*}{Lemma}

\newcommand{\R}{\mathbb{R}}

\DeclarePairedDelimiter{\abs}{\lvert}{\rvert}

\newcommand{\arxiv}[1]{\href{http://arxiv.org/abs/#1}{\texttt{arXiv:#1}}}

\title{A Useful Inequality for the\\ Binary Entropy Function}

\author{Ravi B. Boppana\thanks{Department of Mathematics, Massachusetts Institute of Technology, Cambridge, Massachusetts, USA\@.  Email address: {\tt rboppana@mit.edu}.}
}

\date{January 23, 2023}

\begin{document}

\maketitle

\begin{abstract}
We provide a simple proof of a curious inequality for the binary entropy function,
an inequality that has been used in two different contexts.
In the 1980's, Boppana used this entropy inequality to prove lower bounds on Boolean formulas.
More recently, the inequality was used to achieve major progress on Frankl's union-closed sets conjecture.
Our proof of the entropy inequality uses basic differential calculus.
\end{abstract}

\section{Introduction}

In this note,
we provide a simple proof of a curious inequality for the binary entropy function,
an inequality that has been used in two different contexts.  

Let $h$ be the binary entropy function, 
defined on the interval~$[0, 1]$ as follows:
\[
  h(x) = 
    \begin{cases}
      -x \log x - (1 - x) \log(1 - x) & \text{if  $0 < x < 1$;} \\
      0 & \text{if $x = 0$ or $x = 1$.}
    \end{cases}
\]
Here $\log$ means natural logarithm (base~$e$).

The following lemma is the inequality that we will prove.

\begin{lemma*}
If $0 \le x \le 1$, then
\[
  h(x^2) \ge \phi x h(x) ,
\]
where $\phi$ is the golden ratio, 
$\frac{\sqrt{5} + 1}{2}$.
\end{lemma*}

The history of this entropy inequality is interesting.
Boppana~\cite{Boppana1, Boppana2} 
first used this inequality to prove lower bounds on Boolean formulas.
More precisely,
he proved a two-variable inequality on the binary entropy function,
which as a special case yields the one-variable inequality above.
His proof of the two-variable inequality was computer assisted.
In unpublished work~\cite{Boppana3}, he gave a simple proof of the one-variable inequality.  
This simple proof is the one we will record in this note.

More than 30 years later,
the same inequality was used to make major progress on Frankl's union-closed sets conjecture.  
Gilmer~\cite{Gilmer} used the information-theoretic concept of entropy to achieve a breakthrough on the union-closed conjecture.
Immediately after, 
the entropy inequality above was used to improve Gilmer's bound by
Alweiss, Huang, and Sellke~\cite{AHS},
Chase and Lovett~\cite{CL},
Pebody~\cite{Pebody},
and Sawin~\cite{Sawin}.
Regarding the entropy inequality itself,
Alweiss, Huang, and Sellke~\cite{AHS}
gave a proof using computer assistance.
Chase and Lovett~\cite{CL} cited the proof of~\cite{AHS}.
Pebody~\cite{Pebody} wrote ``to be proven''.  
Sawin~\cite{Sawin} gave a symbolic proof,
noting that it is ``somewhat complicated, though it would not be surprising if a simple symbolic proof exists''.  
Our proof confirms that a simple symbolic proof exists.

Further progress on the union-closed conjecture was given by 
Cambie~\cite{Cambie}, Ellis~\cite{Ellis}, and Yu~\cite{Yu}.

\section{Proof of the entropy inequality}

In this section, 
we provide a proof of the entropy inequality using basic differential calculus.  

\begin{lemma*}
If $0 \le x \le 1$, then
\[
  h(x^2) \ge \phi x h(x) ,
\]
where $\phi$ is the golden ratio, 
$\frac{\sqrt{5} + 1}{2}$.
\end{lemma*}

\begin{proof}
Let $\R$ be the set of real numbers.
It will be convenient to extend~$h$ to all of~$\R$ as follows:
\[
  h(x) = 
    \begin{cases}
      -x \log \abs{x} - (1 - x) \log \abs{1 - x} & \text{if $x \ne 0$ and $x \ne 1$;}  \\
      0 & \text{if $x = 0$ or $x = 1$.}
    \end{cases}
\]

Let $f$ be the function on~$\R$ defined by
\[
  f(x) = h(x^2) - \phi x h(x) .
\] 
To prove the lemma,
we will show that $f$ is nonnegative on~$[0, 1]$.

The first derivative~$f'$ is defined on~$\R \smallsetminus \{ -1, 1 \}$, 
the second derivative~$f''$ is defined on~$\R \smallsetminus \{ -1, 0, 1 \}$,
and the third derivative~$f'''$ is defined on~$\R \smallsetminus \{ -1, 0, 1 \}$.
Taking derivatives three times,
we see that if $x \in \R \smallsetminus \{ -1, 0, 1 \}$, 
then
\[
  f'''(x) = \frac{p(x)}{x (1 - x^2)^2} \, ,
\]
where $p$ is the cubic polynomial defined by
\[
  p(x) = - \phi x^3 - 4 x^2 + 3 \phi x + 2 \phi - 4 .
\]

Because the leading coefficient of~$p$ is negative and $p(0)$ is negative,
$p$ has at least one negative root.
Thus $p$ has at most two nonnegative roots.
Hence the third derivative~$f'''$ has at most two roots in~$(0, 1)$.

By Rolle's theorem,
applied three times,
it follows that $f$ itself has at most five roots in~$[0, 1]$, 
counting multiplicity.
The function~$f$ has a double root at~$0$,
a double root at~$\phi^{-1}$,
and a single root at~$1$.
Thus we have found all five roots of~$f$ in~$[0, 1]$.

Because $f$ has a double root at~$\phi^{-1}$,
it is either all nonnegative or all nonpositive on~$[0, 1]$.  
For $x$ a tiny positive number,
$f(x)$ is positive.
Hence $f$ is nonnegative on~$[0, 1]$.  
\end{proof}

\end{document}